\documentclass{article}

\def\ra{\rightarrow}

\def\ss{\subseteq}

\def\Re{\hbox{\rm Re}\,}

\def\O{{\cal O}}

 \def\HollowBox #1#2{{\dimen0=#1 \advance\dimen0 by -#2       
       \dimen1=#1 \advance\dimen1 by #2                       
        \vrule height #1 depth #2 width #2                    
        \vrule height 0pt depth #2 width #1                   
        \llap{\vrule height #1 depth -\dimen0 width \dimen1}% 
       \hskip -#2                                             
       \vrule height #1 depth #2 width #2}}                   
 \def\BoxOpTwo{\mathord{\HollowBox{6pt}{.4pt}}\;}             
\def\endpf{\hfill $\BoxOpTwo$}

\def\sss{\subset \, \, \subset}

\font\teneufm=eufm10
\font\seveneufm=eufm7
\font\fiveeufm=eufm5
\newfam\eufmfam
\textfont\eufmfam=\teneufm
\scriptfont\eufmfam=\seveneufm
\scriptscriptfont\eufmfam=\fiveeufm

\newfam\msbfam
\font\tenmsb=msbm10  \textfont\msbfam=\tenmsb
\font\sevenmsb=msbm7  \scriptfont\msbfam=\sevenmsb
\font\fivemsb=msbm5    \scriptscriptfont\msbfam=\fivemsb
\def\Bbb{\fam\msbfam \tenmsb}

\def\O{\Omega}

\def\RR{{\Bbb R}}
\def\CC{{\Bbb C}}

\def\ZZ{{\Bbb Z}}

\newtheorem{theorem}{Theorem}[section]

\newtheorem{proposition}[theorem]{Proposition}
\newtheorem{lemma}[theorem]{Lemma}
\newtheorem{remark}[theorem]{Remark}

\makeindex

\usepackage{graphicx}

\usepackage{amsmath}

\begin{document}

\begin{center}
\huge \bf
The Schwarz Lemma at the Boundary
\end{center}
\vspace*{.12in}

\begin{center}
\large Steven G. Krantz\footnote{Author supported in part
by the National Science Foundation and by the Dean of the Graduate
School at Washington University.}\footnote{{\bf MR Classification
Numbers:}  30C80, 30C99, 32A10, 32A30.}
\end{center}
\vspace*{.15in}

\begin{center}
\today
\end{center}
\vspace*{.2in}

\begin{quotation}
{\bf Abstract:} \sl
The most classical version of the Schwarz lemma involves the behavior at the origin
of a bounded, holomorphic function on the disc.  Pick's version of
the Schwarz lemma allows one to move the origin to other points of the disc.

In the present paper we explore versions of the Schwarz lemma at a boundary
point of a domain (not just the disc).   Estimates on derivatives of the function,
and other types of estimates as well, are considered.  We review recent results
of several authors, and present some new theorems as well.
\end{quotation}
\vspace*{.25in}

\setcounter{section}{-1}

\section{Introduction}

The classical Schwarz lemma gives information about the behavior of a holomorphic
function on the disc at the origin, subject only to the relatively mild
hypotheses that the function map the disc to the disc and the origin to the origin.
Later generalizations by Pick allow one to replace ``origin'' by other points
of the disc.  Of course there are far-reaching generalizations of the
classical Schwarz lemma, due to Ahlfors and others, that place the Schwarz lemma
squarely in the province of differential geometry.

In the present paper we explore Schwarz lemmas at the boundary of a domain.  We give
both function-theoretic and geometric formulations of the theorems.  A number of
different proofs and perspectives on the results are presented.

\section{The Classical Schwarz Lemma}

In its most basic form, the familiar Schwarz lemma says this:

\begin{proposition} \sl
Let $f: D \ra D$ be a holomorphic function that fixes the origin 0.
Then 
$$
|f(\zeta)| \leq |\zeta| \qquad \hbox{for all} \ \zeta \in D
$$
and 
$$
|f'(0)| \leq 1 \, .
$$

If $|f(\zeta)| = |\zeta|$ for any $\zeta \ne 0$ or if $|f'(0)| = 1$, then $f$ is
a rotation of the disc.
\end{proposition}
{\bf Proof:}  Apply the maximum principle to the function $g(\zeta) = f(\zeta)/\zeta$.
\endpf
\smallskip \\

Now, as is well known, one may replace the origin in the domain of $f$ and the origin in the range of $f$
in the last proposition to obtain the so-called Schwarz-Pick lemma:

\begin{proposition} \sl
Let $f:D \ra D$ be holomorphic.  Assume that $a \ne b$ are elements of $D$
and that $f(a) = \alpha$, $f(b) = \beta$.  Then
\begin{enumerate}
\item[{\bf (a)}]  $\displaystyle 
 \left | \frac{\beta - \alpha}{1 - \overline{\alpha}\beta} \right |
    \leq \left | \frac{b - a}{1 - \overline{a} b} \right |$;
\item[{\bf (b)}]  $\displaystyle
 \left | f'(a) \right | \leq \frac{1 - |\alpha|^2}{1 - |a|^2}$.
\end{enumerate}
There is also a pair of uniqueness statements:
\begin{enumerate}
\item[{\bf (c)}]  If $\displaystyle 
 \left | \frac{\beta - \alpha}{1 - \overline{\alpha}\beta} \right |
    = \left | \frac{b - a}{1 - \overline{a} b} \right |$,
then $f$ is a conformal self-map of the disk $D$;
\item[{\bf (d)}]  If $\displaystyle
 \left | f'(a) \right | = \frac{1 - |\alpha|^2}{1 - |a|^2}$,
then $f$ is a conformal self-map of the disk $D$.
\end{enumerate}
\end{proposition}

\begin{remark} \rm
The expression
$$
\rho(a,b) = \left | \frac{b - a}{1 - \overline{a} b} \right |
$$
is the {\it pseudohyperbolic metric} on the disc.  Thus {\bf (a)} says that
the mapping $f$ is distance decreasing in the pseudohyperbolic metric.
It is noteworthy that the pseudohyperbolic metric is {\it not} a Riemannian
metric.
\end{remark}

\noindent {\bf Proof of the Proposition:}   Recall that, for
$a$ a complex number in $D$,
$$
\varphi_a(\zeta) = \frac{\zeta - a}{1 - \overline{a}\zeta}
$$
defines a {\small \it M\"{o}bius transformation}.  This is a conformal
self-map of the disk that takes $a$ to 0.  Note that $\varphi_{-a}$
is the inverse mapping to $\varphi_a$.

Now, for the given $f$, consider
$$
g(\zeta) = \varphi_\alpha \circ f \circ \varphi_{-a} (\zeta) \, .
$$
Then $g: D \ra D$ and $g(0) = 0$.  So the standard Schwarz
lemma applies to $g$.  By part {\bf (a)} of that lemma,
$$
|g(\zeta)| \leq |\zeta| \, .
$$
Letting $\zeta = \varphi_a(\xi)$ yields
$$
|\varphi_\alpha \circ f(\xi)| \leq |\varphi_a(\xi)| \, .
$$
Writing this out, and setting $\zeta = b$, gives the conclusion
$$
\left | \frac{\beta - \alpha}{1 - \overline{\alpha}\beta} \right |
    \leq \left | \frac{b - a}{1 - \overline{a} b} \right | \, .
$$
That is part {\bf (a)}.

For part {\bf (b)}, we certainly have that
$$
\left | (\varphi_\alpha \circ f \circ \varphi_{-a} )'(0) \right | \leq 1 \, .
$$
Using the chain rule, we may rewrite this as
$$
\left | \varphi'_\alpha(f \circ \varphi_{-a}(0)) \right | \cdot
   \left | f'(\varphi_{-a}(0)) \right | \cdot
  \left | \varphi'_{-a} (0) \right | \leq 1 \, .  \eqno (1.2.1)
$$
Now of course 
$$
\varphi'_a(\zeta) = \frac{1 - |a|^2}{(1 - \overline{a}\zeta)^2} \, .
$$
So we may rewrite (1.2.1) as
$$
\left ( \frac{1 - |\alpha|^2}{(1 - |\alpha|^2)^2} \right ) \cdot 
     |f'(a)| \cdot (1 - |a|^2) \leq 1 \, .
$$
Now part {\bf (b)} follows.

We leave parts {\bf (c)} and {\bf (d)} as exercises for the reader.
\endpf 
\smallskip \\

It is easy to see that the statement of the Schwarz-Pick lemma degenerates as the point
$a$ tends to the boundary.  So some other idea will be required if we are to
successfully formulate and prove a boundary Schwarz lemma.

\section{A First Look at the Boundary Schwarz Lemma}

\subsection{The Hopf Lemma}

The next result is one of the antecedents to a classical Schwarz lemma at the boundary.  We shall
first state the lemma, then discuss its context and significance.

\begin{lemma}[Hopf]  \sl
Let $\Omega \sss \RR^N$ have $C^2$ boundary.  Let $u \in C(\overline{\Omega})$ be real-valued
with $u$ harmonic and non-constant on $\Omega.$  Let $P \in \partial \Omega$ and assume that
$u$ takes a local minimum at $P.$  Then
$$
\frac{\partial u}{\partial \nu} (P) < 0 . 
$$
\end{lemma}
{\bf Proof:}  Suppose without loss of generality that
$u > 0$ on $\Omega$ near $P$ and that 
$u(P) = 0.$  Let $B_R$ be a ball that is internally
tangent to $\partial \Omega$ at $P.$  We may assume that
the center of this ball is at the origin and that
$P$ has coordinates $(R,0,\dots,0).$
Then, by Harnack's inequality
(see [KR1]), we have for $0 < r < R$ that
$$
u(r,0,\dots,0) \geq c \cdot \frac{R^2 - r^2}{R^2 + r^2}
$$
hence
$$
\frac{u(r,0,\dots,0) - u(R,0,\dots,0)}{r - R} \leq - c' < 0 .
$$
Therefore
$$
\frac{\partial u}{\partial \nu}(P) \leq - c' < 0 . 
$$
This is the desired result.
\endpf 
\smallskip \\

A good reference for the Hopf lemma is [COH].  It was used in
that source to provide a proof of the maximum principal for second-order,
elliptic partial differential operators.  Namely, if a solution	$u$
of such an operator ${\cal L}$ has an interior maximum at a point $P$, then
let $S$ be a sphere passing through $P$.  Restrict attention to the closed
ball $B$ bounded by $S$.   Then the function $u$ has a maximum at $P$, so the
outward normal derivative at $P$ is positive.  But that means that, at a point
near $P$ in the outward normal direction the function $u$ will take an even larger value,
contradicting the maximality of $u$ at $P$.

In more recent times the Hopf lemma has proved particularly useful in the study
of biholomorphic and proper holomorphic mappings of several complex variables (see, for instance, [KRA1]).

The Hopf lemma is true in fact for subharmonic functions, and under rather weak hypotheses
on the behavior of $u$ at $P$.  We leave the details for the interested reader.  The message
that the Hopf lemma gives us is best seen for a holomorphic mapping $F: B \ra B$, where
$B$ is the unit ball in $\CC^n$.  Let ${\bf 1} = (1,0,0, \dots, 0) \in \partial B$, and assume
that the limit of $F(z)$ is ${\bf 1}$ as $z$ approaches {\bf 1} admissibly (see [KRA1] for this
concept).   Let $\nu$ be the unit normal vector to the boundary at {\bf 1}, and set
$f(z) = F(z) \cdot \nu$.  Finally let $u(z) = |f(z)|$.  Then $u$ is plurisubharmonic, and
$u$ takes a maximum value (in a reasonable sense) at {\bf 1}.  The Hopf lemma applies, and we
see that the normal derivative of $u$ at {\bf 1} is nonzero.  This tells us that the
boundary point {\bf 1} is analytically isolated for the function $f$.  And that is
a primitive version of the Schwarz lemma at the boundary point {\bf 1}.

\section{Work of L\"{o}wner and Velling}

As early as 1923, K. L\"{o}wner was considering deformation theorems that can be considered
to be early versions of the Schwarz lemma at the boundary.   A version of his result is this:

\begin{proposition} \sl Let $f: D \ra D$ holomorphic with $f(0)
= 0$. Of course $f$ has radial boundary limits almost
everywhere. Let $S = \partial D$. Assume that $f$ maps an arc
$A \ss S$ of length $s$ onto an arc $f(A) \ss S$ of length
$\sigma$. Then $\sigma \geq s$ with equality if and only if
either $s = \sigma = 0$ or $f$ is just a rotation.
\end{proposition}

We see that L\"{o}wner tells us that a boundary arc must be mapped to a boundary arc that is shorter.
This is in the spirit of the original Schwarz lemma, for it tells us that, under a similar mapping,
the stretching factor must be less than 1.

John Velling studied L\"{o}wner's ideas in 1985 and proved a refinement which we shall treat
at the end of the next section.

\section{A Refinement}

In [OSS], R. Osserman offered the following boundary refinement of the classical Schwarz lemma.
It is very much in the spirit of the sort of result that we wish to consider here.

\begin{theorem}  \sl
Let $f: D \ra D$ be holomorphic.  Assume that $f(0) = 0$.  Further assume that
there is a $b \in \partial D$ so that $f$ extends continuously to $b$,
$|f(b)| = 1$ (say that $f(b) = c$), and $f'(b)$ exists.  Then
$$
|f'(b)| \geq \frac{2}{1 + |f'(0)|} \, .	   \eqno (4.1.1)
$$
\end{theorem}

\begin{remark} \rm 
It is easy to see that inequality $(4.1.1)$ is sharp.  For instance, the
function $f(\zeta) =\zeta$ gives equality.  In fact, for each possible
value of $|f'(0|$ (between 0 and 1 inclusive) there is a function
that makes $(4.1.1)$ sharp.   For $0 \leq a \leq 1$, the
function 
$$
f(\zeta) = \zeta \cdot \frac{\zeta + a}{1 + a \zeta}
$$
gives equality in $(4.1.1)$.
\end{remark}

By way of proving the theorem, we first prove some preliminary results.

\begin{lemma} \sl
Let $f: D \ra D$ be holomorphic and satisfy $f(0) = 0$.  Then
$$
|f(\zeta)| \leq |\zeta| \cdot \frac{|\zeta| + |f'(0)|}{1 + |f'(0)| |\zeta|} \qquad \hbox{for} \ |\zeta| < 1 \, .  \eqno (4.3.1)
$$
\end{lemma}
{\bf Proof:}  As usual, set $g(\zeta) = f(\zeta)/\zeta$.  The usual Schwarz lemma then tells
us that either $f$ is a rotation or else $|g(\zeta)| < 1$ for $|\zeta| < 1$.  The first
of these eventualities leads to $|f'(0)| = 1$ and hence our inequality is trivially true.
So we may as well suppose that $|g(\zeta)| < 1$ for $|\zeta| < 1$.   Applying a rotation if
necessary, we may also suppose that $g(0) = f'(0) = a$, where $a$ is real and $0 \leq a < 1$.

Thus inequality $(4.3.1)$ is equivalent to
$$
|g(\zeta)| \leq \frac{|\zeta| + a}{1 + a |\zeta|} \qquad \hbox{for} \ |\zeta| < 1 \, .  \eqno (4.3.2)
$$
We may derive this assertion from the Schwarz-Pick lemma as follows.  Certainly $g$ will
map each disc $D(0,r)$, $0 < r < 1$, into a disc with diameter the real interval
$$
\left [ \frac{a - r}{1 - ar}, \frac{a + r}{1 + ar} \right ] \, .
$$
As a result, when $|\zeta| = r$ then 
$$
|g(\zeta)| \leq \frac{a + r}{1 + ar} = \frac{|\zeta| + a}{1 + a|\zeta|} \, .
$$
That proves $(4.3.2)$.  Then $(4.3.1)$ follows.
\endpf
\smallskip \\

\begin{remark} \rm
In view of the second part of the classical Schwarz lemma, the fraction
$[|\zeta| + |f'(0)|]/[1 + |f'(0)||\zeta|]$ does not exceed 1.  Thus
one recovers the basic inequality of the usual Schwarz lemma.
\end{remark}

\begin{lemma} \sl
We have
$$
\lim_{\zeta_j \ra b} \left | \frac{f(\zeta_j) - c}{|\zeta_j| - |b|} \right | \geq
    \lim_{\zeta_j \ra b} \frac{1 - |f(\zeta_j)|}{1 - |\zeta_j|} \geq \frac{2}{1 + |f'(0)|} \, .
$$
\end{lemma}
{\bf Proof:}  Certainly
$$
\left | \frac{f(\zeta) - c}{|\zeta| - |b|} \right | \geq \frac{1 - |f(\zeta)|}{1 - |\zeta|}
$$
by elementary inequalities (using of course the facts that $|b| = 1$ and $|c| = 1$).  But now
the last lemma gives an upper bound for $|f(\zeta)|$.  Plugging that into
the righthand side of this last inequality gives the second inequality.
\endpf
\smallskip \\

\noindent {\bf Proof of the Theorem:}  Simply let $\zeta_j$ in the last lemma
equal $t_j  b$ for real $t_j \ra 1$.  Letting $j \ra +\infty$ then gives that the
lefthand side becomes $|f'(b))|$ and the
result follows.  
\endpf
\smallskip \\

A very interesting consequence of the theorem is the following result, which is a refinement of a theorem of Velling [VEL].

\begin{proposition} \sl
Let $f: D \ra D$ be holomorphic.  Let $S \ss \partial D$ be a nontrivial arc, and suppose that $f$ extends
continuously to $S$.  Further assume that $f(S)$ lies in $\partial D$. 	 Let $s$ denote the length of $S$ and
$\sigma$ the length of $f(S)$ (which is also necessarily an arc, since it is a connected subset of the circle).  Then
$$
\sigma \geq \frac{2}{1 + |f'(0)|} \cdot s \, .
$$
\end{proposition}
{\bf Proof:}  By Schwarz reflection, we may take it that $f$ is analytic on the interior of the arc $S$.  Hence it certainly
satisfies the hypotheses of the first lemma at each point of the interior of $S$.    The conclusion of that lemma then holds,
and integration yields the desired result.
\endpf
\smallskip \\

\section{A New Look at the Schwarz Lemma on the Boundary}

Certainly the most interesting and valuable part of the Schwarz lemma is the uniqueness statement for the derivative.  {\it That}
is the result for which we would like to have a boundary formulation.  The next theorem is due to Burns and Krantz [BUK].

\begin{theorem}[Burns/Krantz]
Let $\phi: D \ra D$ be a holomorphic function from the disc to itself such that
$$
\phi(\zeta) = 1 + (\zeta - 1) + {\cal O}\bigl (|\zeta - 1|^4) \bigr )  
$$
as $\zeta \ra 1$.  Then $\phi(\zeta) \equiv \zeta$ on the disc.
\end{theorem}

\begin{remark} \rm
Today there are several proofs of this result.   Chelst [CHE] has some nice ways to look at the
matter.   Boas [BOA] has a new and brief proof.  We present here the original proof because
it is somewhat natural and also enlightening.

It is worth noting that Velling [VEL1] and others have proved antecedents to
this theorem.  But their results had additional hypotheses, such as univalence
of the function, or analyticity in a neighborhood of 1.   The result presented
here is the sharpest possible.  Indeed, the example
$$
\phi(\zeta) = \zeta - \frac{1}{10} \cdot (\zeta - 1)^3
$$
shows that the exponent 4 in the theorem cannot be replaced by 3.   The proof
in fact shows that 4 {\it can} be replaced by $o\bigl ( |\zeta - 1|^3 \bigr )$.
\end{remark}

\noindent {\bf Proof of the Theorem:}   Consider the holomorphic function
$$
g(\zeta) = \frac{1 + \phi(\zeta)}{1 - \phi(\zeta)} \, .
$$
Then $g$ maps the disc $D$ to the right halfplane.  By the Herglotz representation (see [AHL2]), there must be a positive 
measure $\mu$ on the interval $[0, 2\pi)$ and an imaginary constant ${\cal C}$ so that
$$
g(\zeta) = \frac{1}{2\pi} \int_0^{2\pi} \frac{e^{i\theta} + \zeta}{e^{i\theta} - \zeta} \, d\mu(\theta) + {\cal C} \, .  \eqno (5.1.1)
$$

The hypothesis on $\phi$ will enable us to analyze the structure of $g$ and hence the structure of $\mu$.  
To wit, we write
$$
g(\zeta) = \frac{1 + \zeta + {\cal O}(\zeta - 1)^4}{1 - \zeta - {\cal O}(\zeta - 1)^4} = \frac{1 + \zeta}{1 - \zeta} + {\cal O}(\zeta - 1)^2 \, .
$$
This and equation $(5.1.1)$ imply that the measure $\mu$ has the form $\mu = \delta_0 + \nu$, where $\delta_0$ is ($2\pi$ times) the Dirac mass
at the origin and $\nu$ is another positive measure on $[0, 2\pi)$.  In fact a good way to verify the positivity of $\nu$ is to use the
equation
$$
\frac{1 + \zeta}{1 - \zeta} + {\cal O}(\zeta - 1)^2 = \frac{1}{2\pi} \int_0^{2\pi}  \frac{e^{i\theta} + \zeta}{e^{i\theta} - \zeta} \, d(\delta_0 + \nu) (\theta) + {\cal C} 
$$
to derive a Fourier-Stieltjes expansion of $\delta_0 + \nu$ and then to apply the Herglotz criterion [KAT, p.\ 38].

We may simplify this last equation to
$$
{\cal O}(\zeta - 1)^2 = \frac{1}{2\pi} \int_0^{2\pi} \frac{e^{i\theta} + \zeta}{e^{i\theta} - \zeta} \, d\nu(\theta) + {\cal C} \, .
$$
Now pass to the real part of the last equation.  That eliminates the constant ${\cal C}$.  Since $\nu$ is a positive measure, we see that the real
part of the integral on the righthand side of this last equation represents a positive harmonic function $h$ on the disc that
satisfies
$$
h(\zeta) = {\cal O}(\zeta - 1)^2 \, .
$$
In particular, $h$ takes a minimum at the point $\zeta = 1$ and is ${\cal O}(|\zeta - 1|^2)$.  This contradicts
Hopf's lemma (see [KRA1, [GRK]) unless $h \equiv 0$.  But $h \equiv 0$ means that $\nu \equiv 0$.  Therefore
$$
g(\zeta) = \frac{1 + \zeta}{1 - \zeta} \, .
$$
We conclude that $\phi(\zeta) \equiv \zeta$.  That is the assertion that we wish to prove.
\endpf
\smallskip \\

It is worthwhile to formulate the classical Schwarz lemma in the language
of this last theorem.  One way to do this is as follows:
\begin{lemma} \sl
Let $f: D \ra D$ be holomorphic, and assume that $f(0) = 0$.  If
$$
f(\zeta) = \zeta + {\cal O}(|\zeta|^2) \, ,    \eqno (5.3.1)
$$
then $f(\zeta) \equiv \zeta$.   
\end{lemma}

\noindent The proof is obvious, for the hypothesis $(5.3.1)$ implies
that $f'(0) = 1$.   

We might also recall H. Cartan's classic result:

\begin{theorem} \sl
Let $\Omega \ss \CC^n$ be a bounded domain.  Fix a point $P \in \Omega$.  Suppose
that $\phi : \Omega \ra \Omega$ is a holomorphic mapping such that $\phi(P) = P$.  If 
the complex Jacobian of $\phi$ at $P$ is the identity matrix, then $\phi$ is the identity
mapping.
\end{theorem}

\noindent We may think of Cartan's theorem as a reformulatin of (5.3.1) in the multivariable setting.
We now, for the sake of interest and completeness, provide a proof of Cartan's result.
\smallskip \\

\noindent {\bf Proof of Theorem 5.4:}   We may assume that $P = 0.$  Expanding $\phi$ in 
a power series about $P = 0$ (and remembering that $\phi$ is vector-valued
hence so is the expansion) we have
$$
  \phi(z) = z + P_k(z) + O(|z|^{k+1}) , 
$$
where $P_k$is the first homogeneous polynomial of order exceeding
 $1$ in the Taylor expansion.  Defining $\phi^j(z) = \phi \circ \cdots \circ \phi$
($j$ times) we have
\begin{eqnarray*}
  \phi^2(z) & = & z + 2P_k(z) + O(|z|^{k+1}) \\
  \phi^3(z) & = & z + 3P_k(z) + O(|z|^{k+1}) \\
            & \cdot &                            \\
            & \cdot &                            \\
            & \cdot &                            \\
  \phi^j(z) &   =   & z + jP_k(z) + O(|z|^{k+1}) .
\end{eqnarray*}

Choose polydiscs $D^n(0,a) \ss \O \ss D^n(0,b).$  Then
for $0 \leq j \in \ZZ$ we know that $D^n(0,a) \ss \mbox{dom}\, \phi^j \ss D^n(0,b).$
Therefore the Cauchy estimates imply that for any multi-index $\alpha$ with
$|\alpha| = k$ we have
$$
 j|D^\alpha \phi(0)| = |D^\alpha \phi^j(0)| \leq n \frac{b \cdot \alpha!}{a^k} . 
$$
Letting $j \ra \infty$ yields that $D^\alpha \phi(0) = 0.$

We conclude that $P_k = 0;$ this contradicts the choice of $P_k$ unless
$\phi(z) \equiv z.$ 
\endpf
\smallskip \\

\begin{remark} \rm
Notice that this proposition is a generalization
of the uniqueness part of the classical Schwarz lemma on the disc.  In fact
a great deal of work has been devoted to generalizations of this type of Schwarz
lemma to more general settings.  We refer the reader to 
[WU], [YAU], [KRA4], [KRA5], [BUK] for more on this matter.
\end{remark}

\section{Ideas of Chelst}

The following lemma is relevant to our considerations in this section.

\begin{lemma} \sl
Let $\Omega$ be a bounded domain in $\CC$ and let $u$ be a real-valued
harmonic function on $\Omega$.  Suppose that there is a collared neighborhood
$U$ of $\partial \Omega$ so the $u \geq 0$ on $U \cap \Omega$.  Then $u \geq 0$
everywhere.
\end{lemma}

\begin{remark} \rm
It is not enough for $u$ to simply be nonnegative on $\partial \Omega$.  As a simple
example, let $\Omega$ be the upper halfplane and let $u(x,y) = x^2 - y^2$.  Then
clearly $u \geq 0$ on $\partial \Omega$---indeed $u > 0$ at every point of $\partial \Omega$ except
the origin.  Yet $u$ is {\it not} nonnegative on the positive imaginary axis.
\end{remark}

\noindent {\bf Proof of the Lemma:}  Applying the maximum principle to $-u$ on a slightly
smaller domain (with boundary lying inside $U \cap \Omega$), we see that $-u$ cannot be
positive in $\Omega \setminus U$.  Hence $u \geq 0$ on all of $\Omega$.
\endpf 
\smallskip \\

We will also make good use of the classical Hopf lemma, as enunciated in an earlier
part of the present paper.  

Now the following proposition is inspired by Chelst's main result, but is strictly more general.
As a result, the line of argument is necessarily different.

\begin{proposition} \sl  
Let $f:D \ra D$ be a holomorphic function.  Let $B$ be an inner function which
equals 1 precisely on a set $A_B \ss \partial D$ of measure 0.  Assume that
\begin{enumerate}
\item[{\bf (a)}]  For a given point $a \in A_B$, $f(\zeta) = B(\zeta) + {\cal O}(|\zeta - a|^4)$ as $\zeta \ra a$;
\item[{\bf (b)}]  For all $b \in A_B \setminus \{a\}$, $f(\zeta) = B(\zeta) + {\cal O}(|\zeta - b|^2)$ as $\zeta \ra b$.
\end{enumerate}
Then $f(\zeta) \equiv B(\zeta)$ on all of $D$.
\end{proposition}

\begin{remark} \rm
It needs to be clearly understood here that $A_B$ is the full set on which $B$ equals 1.  The proof
consists of coming to terms with the boundary behavior of $f$ and $B$ on that set.
\end{remark}

{\bf Proof:}  	Following Chelst, it is useful to consider the function
$$
h(\zeta) = \Re \left [ \frac{1 + f(\zeta)}{1 - f(\zeta)} \right ] - \Re \left [ \frac{1 + B(\zeta)}{1 - B(\zeta)} \right ]  \, .
$$
We shall perform some estimates to show that {\bf (i)}  $h$ has non-negative boundary limits almost everywhere
on $\partial D$ and {\bf (ii)}	$h$ lies in $h^2(D)$ (i.e., harmonic functions which are uniformly
square integrable on circles centered at the origin---see [KRA1]).  The natural conclusion then
is that $h$ is positive everywhere on the interior of $D$.

Now
\begin{eqnarray*}
h(\zeta) & = & \Re \left [ \frac{1 + f(\zeta)}{1 - f(\zeta)} \right ] - \Re \left [ \frac{1 + B(\zeta)}{1 - B(\zeta)} \right ]	   \\
	 & = & \Re \left [ \frac{[1 + B(\zeta) + {\cal O}(|\zeta - 1|^4)] \cdot [1 - \overline{B(\zeta)} + {\cal O}(|\zeta - 1|^4)]}{|1 - B(\zeta) + {\cal O}(|\zeta - 1|^4)|^2} \right ] 
		  - \Re \left [ \frac{(1 + B(\zeta))(1 - \overline{B(\zeta)})}{|1 - B(\zeta)|^2}  \right  ]  \\
	 & = & \Re \left [ \frac{(1 - \overline{B(\zeta)} + B(\zeta) - |B(\zeta)|^2 +  {\cal O}(|\zeta - 1|^4)}{|1 - B(\zeta) + {\cal O}(|\zeta - 1|^4)|^2} \right ] 
                  - \Re \left [ \frac{(1 - \overline{B(\zeta)} + B(\zeta) - |B(\zeta)|^2}{|1 - B(\zeta)|^2} \right ]     \\
	 & = & \frac{[1 - |B(\zeta)|^2 + {\cal O}(|\zeta - 1|^4)] \cdot |1 - B(\zeta)|^2 - [|1 - B(\zeta) + {\cal O}(|\zeta - 1|^4) ]^2 \cdot ( 1 - |B(\zeta)|^2)}{|1 - B(\zeta) + {\cal O}(|\zeta - 1|^4)|^2 \cdot |1 - B(\zeta)|^2} \\
	 & = & \frac{[(1 - |B(\zeta)|^2) \cdot |1 - B(\zeta)|^2 + {\cal O}(|\zeta - 1|^4 ] - [|1 - B(\zeta)|^2 \cdot (1 - |B(\zeta)|^2) + {\cal O}(|\zeta - 1|^4)]}{|1 - B(\zeta)|^4} \\
         & = & \frac{{\cal O}(|\zeta - 1|^4)}{|1 - B(\zeta)|^4}  \, . \\
\end{eqnarray*}
But Hopf's lemma tells us that $|1 - B(\zeta)|$ is {\it not} $o(|\zeta - 1|)$.  And in fact we can certainly
say (a bit sloppily) that $|1 - B(\zeta)| \geq C \cdot |1 - \zeta|^{1 + \epsilon}$ for some small $\epsilon > 0$.

In conclusion, the function $h$ certainly lies in $h^2(D)$.   We also note that (and our calculations show this)
the boundary limits of the first expression on the righthand side of the first line of the previous multi-line display are nonnegative almost everywhere.
And the boundary limits of the second expression on the righthand side of the first line of the previous multi-line display are 0 almost everywhere.
In summary, we have an $h^2$ harmonic function with nonnegative radial boundary limits almost everywhere.  It then follows, from
the Poisson integral formula for instance, that $h$ is positive on the disc $D$.  

But $h$ takes the boundary limit 0 at each point of $A_B$.  It follows then from Hopf's lemma that $h$ has a nonzero normal
derivative at each of those points.  That fact contradicts hypothesis {\bf (b)} of the proposition.  And that contradiction
tells us that $h \equiv 0$ hence $f$ is identically equal to the Blaschke product $B$.
\endpf 
\smallskip \\

Chelst [CHE] has pointed out that the function
$$
f(\zeta) = \zeta^8 - \frac{1}{256} (\zeta + 1) \bigl [ (\zeta^2 + 1) (\zeta^4 + 1) \bigr ]^2 \cdot (\zeta - 1)^4
$$
maps the disc to the disc and fails hypothesis {\bf (b)} of Proposition 6.3 with $A_B = \{-1,1\}$ and
$B = \zeta \cdot \zeta$; it also fails the conclusion.

It should be mentioned that the papers [BZZ] and [SHO] offer further refinements of the Burns/Krantz and Chelst
theorems.

\section{Variants in the Several Complex Variables Setting}

The work described above, in the one-complex-variable setting, from [BUK] was inspired by a question
of several complex variables.  Namely one wanted to know whether a holomorphic mapping $\Phi: B \ra B$ (where $B$ is
the unit ball in $\CC^n$) could have boundary image $\Phi(\partial B)$ with high order of contact with the target boundary $\partial B$.
In one complex variable, the Riemann mapping theorem tells us that, for a holomorphic mapping $\varphi: D \ra D$,
any order of contact of $\varphi(\partial D)$ with the target boundary $\partial D$ is possible.  Of course there is
no Riemann mapping theorem in several complex variables, and this together with other {\it ad hoc} evidence
suggested that there ought to be an upper bound on the order of contact in the multi-dimensional case.  

The first step in understanding this situation is to prove a multi-dimensional version of Theorem 5.1:

\begin{theorem} \sl
Let $\Phi: B \ra B$ be a holomorphic mapping.  Let ${\bf 1}  \equiv (1,0,0, \dots, 0)$ be the usual boundary
point of the ball.  Assume that
$$
\Phi(z) = {\bf 1} + (z - {\bf 1}) + {\cal O} (|z - {\bf 1}|^4) \, .
$$
Then $\Phi(z) \equiv z$ for all $z \in B$.
\end{theorem}
{\bf Proof:}    For simplicity we restrict attention to complex dimension 2.  For each $a \in B$, let ${\cal L}_a$ be the complex line passing through $a$ and ${\bf 1}$.  
Let ${\bf d}_a$ be the complex disc given by ${\cal L}_a \cap B$.  With $a$ fixed, consider the
holomorphic function
\begin{eqnarray*}
\psi: D & \longrightarrow & B  \\
      \zeta & \longmapsto & (\zeta, 0) \, .
\end{eqnarray*}
Also consider the mapping
$$
\phi_a: B \ra B
$$
which is the automorphism of the ball $B$ which maps ${\bf d}_0$ onto ${\bf d}_a$ and fixes ${\bf 1}$.
Indeed one may say rather explicitly what this last automorphism is.  Note that,for $\alpha$ a complex number
of modulus less than 1, the mapping
$$
\lambda_\alpha(z_1, z_2) = \left ( \begin{array}{c}
                            \displaystyle \frac{(1 - |\alpha|^2)z_1}{1 + \overline{\alpha}z_2} + \frac{\overline{\alpha} (z_2 + \alpha)}{1 + \overline{\alpha} z_2} \\  [.2in]
			    \displaystyle \frac{- \alpha \sqrt{1 - |\alpha|^2} z_1}{1 + \overline{\alpha}z_2} + \frac{(z_2 + \alpha)\sqrt{1 - |\alpha|^2}}{1 + \overline{\alpha} z_2} 
		      	      \end{array}
		      \right )
$$
sends the complex line ${\bf d}_0$ through $(0,0)$ and $(1,0)$ to the complex line through $(|\alpha|^2, \alpha\sqrt{1 - |\alpha|^2})$.   Composition with unitary mappings
will allow us to replace $(|\alpha|^2, \alpha\sqrt{1 - |\alpha|^2})$ with any other element of the ball $B$.

Finally define
\begin{eqnarray*}
\pi_1: B & \longrightarrow & B	\\
       (z_1, z_2) & \longmapsto & (z_1, 0)
\end{eqnarray*}
and
\begin{eqnarray*}
\eta : {\bf d}_0 & \longrightarrow & D  \\
	 (z_1, 0) & \longmapsto & z_1 \, .
\end{eqnarray*}

The function
\begin{eqnarray*}
H_a : D & \longrightarrow & D \\
     \zeta & \longmapsto & \eta \circ \pi_1 \circ ( \phi_a)^{-1} \circ \Phi \circ \phi_a \circ \psi(\zeta)
\end{eqnarray*}
is well defined.    In addition, $H$ satisfies the hypotheses of Theorem 5.1.  It follows then that $H_a(\zeta) \equiv \zeta$.

Now set
$$
G_a(\zeta) = (\phi_a)^{-1} \circ \Phi \circ \phi_a \circ \psi(\zeta) \equiv \bigl ( g_a^1(\zeta), g_a^2(\zeta) \bigr ) \, .
$$
The statement that $H_a(\zeta) \equiv \zeta$ tells us that $g_a^1(\zeta) \equiv \zeta$.  But then
$$
|g_a^1(\zeta)|^2 + |g_a^2(\zeta)|^2 < 1
$$
for $\zeta \in D$.    

Letting $|\zeta| \ra 1$ now yields that $|g_a^2(\zeta)| \ra 0$.  Thus $g_a^2 \equiv 0$.  It now follows
that the image of $G_a$ already lies in ${\bf d}_0$.  Consequently it must be that $\Phi$ preserves ${\bf d}_a$.  This
last assertion can hold for every choice of $a$ if an only if $\Phi$ is the identity mapping.

That completes the proof.
\endpf 
\smallskip \\

It is naturally desirable to extend this last result to a more general class of domains.  The key insight
here is to note that the discs ${\bf d}_a$ in $B$ may be replaced, in a more general setting, by
extremal discs for the Kobayashi metric (see, for instance [KRA1] and especially [LEM]).    The theory of such discs
is well developed in the context of strongly convex domains, and the proof we have given 
here transfers naturally to that setting.

For strongly pseudoconvex domains, there is no theory of extremal discs in the sense of Lempert (but see [KRA3]).  However, Burns
and Krantz [BUK] were able to construct a local theory of extremal discs near a strongly pseudoconvex boundary
point.  As a result, it is possible to prove a version of Theorem 7.1 on a smoothly bounded, strongly pseudoconvex
domain.  Details may be found in [BUK].  For the record, we record the result now:

\begin{theorem} \sl
Let $\Omega \ss \CC^n$ be a smoothly bounded, strongly pseudoconvex domain.
Let $\Phi: \Omega \ra \Omega$ be a holomorphic mapping.  Let $P \in \partial \Omega$ be a boundary
point.  Assume that
$$
\Phi(z) = P + (z - P) + {\cal O} (|z - P|^4) \, .
$$
Then $\Phi(z) \equiv z$ for all $z \in \Omega$.
\end{theorem}

We close this section by using Proposition 6.3 to derive a new version of Theorem 7.1.

\begin{theorem} \sl Let $f:B \ra B$ be a holomorphic function.
Let $h$ be an inner function which equals 1 on a set $A_h \ss
\partial B$ of measure 0. Assume that 
\begin{enumerate}
\item[{\bf (a)}] For a given point $a \in A_h$, $f(z) =
B(z) + {\cal O}(|z - a|^4)$ as $z \ra a$;
\item[{\bf (b)}] For all $b \in A_h \setminus \{a\}$,
$f(z) = B(z) + {\cal O}(|z - b|^2)$ as $z \ra
b$. 
\end{enumerate} Then $f(z) \equiv h(z)$ on all of
$D$. \end{theorem}
{\bf Proof:}  This result is derived from Proposition 6.3 in just the same way that Theorem 7.1 is derived
from Theorem 5.1.
\endpf
\smallskip \\

It is worth mentioning that the work in [FEF] shows that the hypothesis of Theorem 7.2 implies
that the bounday asymptotics of the Bergman metric are preserved (asymptotically at $P$) by
the mapping $\Phi$.  In particular, pseudo-transversal geodesics (in the language of Fefferman)
are mapped to pseudo-transversal geodesics.  And the asymptotic expansion for the Bergman
kernel is mapped to itself in a natural way.

\section{Non-Equidimensional Mappings}

In view of recent work by Webster [WEB], Cima and Suffridge [CIS1], [CIS2],
D'Angelo [DANG1], [DANG2] and others, it is natural to ask what results
may be obtained for mappings $\Phi: \Omega_1 \ra \Omega_2$ where $\Omega_1
\ss \CC^n$, $\Omega_2 \ss \CC^m$, and $n < m$.  In this circumstance
the Levi form, and particularly the type (in the sense of Kohn/D'Angelo/Catlin---see [KRA1]),
is the determining factor.

We begin with a basic result:

\begin{proposition} \sl
Let $\Omega \ss \CC^m$ be a  smoothly bounded domain with defining function $\rho$.  Further let
$\varphi: D \ra \Omega$ be a holomorphic mapping.  Let $P \in \partial \Omega$
be a strongly pseudoconvex point and suppose that
$$
\rho(\varphi(\zeta)) = o (\|\varphi(\zeta) - P\|^2) 
$$
as $\zeta \ra 1 \in \partial D$.  Then $\varphi(\zeta) \equiv P$.
\end{proposition}
{\bf Proof:}  This is simply a restatement of a well-known fact
about a boundary point of type 2 (again see [KRA1] for the definition of,
and discussion of, type).   A strongly pseudoconvex point $P$ is of type
2, hence cannot have a nontrivial analytic disc with order of contact to the 
boundary at $P$ exceeding 2.
\endpf 
\smallskip \\

\begin{proposition} \sl
Let $\Omega \ss \CC^2$ be a  smoothly bounded domain with defining function $\rho$.  Further let
$\varphi: D \ra \Omega$ be a holomorphic mapping.  Let $P \in \partial \Omega$
be a point of geometric type $m$ (see [KRA1, p.\ 468]) and suppose that
$$
\rho(\varphi(\zeta)) = o (\|\varphi(\zeta) - P\|^m) 
$$
as $\zeta \ra 1 \in \partial D$.  Then $\varphi(\zeta) \equiv P$.
\end{proposition}
{\bf Proof:}  The argument is the same as for the last proposition.
\endpf 
\smallskip \\

There are analogous results in higher dimensions, but they are more difficult
to formulate because the concept of type (due to D'Angelo [DANG3]) is more subtle.
We leave the details for another time.

\section{Further Generalizations}

In the paper [HUA], X. Huang was able to generalize the Burns/Krantz
theorem 7.1 to a class of weakly pseudoconvex domains. His key idea---one
that will no doubt see good use in the future---is to use the exponent of
the bounded, plurisubharmonic exhaustion function of [DIF] as a measure of
the geometry of the boundary point. It remains to be seen what the optimal
version of Theorem 7.1 will be on any smoothly bounded pseudoconvex domain
in $\CC^n$.

\section{Concluding Remarks}

The idea of Schwarz lemmas at the boundary has seen considerable activity
in the past ten years or so. It is clearly a providential course of
inquiry, and important for geometric function theory. And there is much
yet to be known. We hope that this paper will point in some new
directions, and inspire some new results.

\newpage
	     
\noindent {\Large \sc References}
\bigskip  \\

\begin{enumerate} 

\item[{\bf [AHL1]}]  L. Ahlfors, An extension of Schwarz's
lemma {\it Trans. Amer. Math. Soc.} 43(1938), 359--364.
				
\item[{\bf [AHL2]}]  L. Ahlfors, {\it Conformal Invariants}, McGraw-Hill, New York, 1973.

\item[{\bf [BZZ]}] L. Baracco, D. Zaitsev, G. Zampieri, A
Burns-Krantz type theorem for domains with corners, {\it
Math.\ Ann.} 336(2006), 491--504.

\item[{\bf [BOA]}]  H. Boas, Julius and Julia:  Mastering the art of
the Schwarz lemma, preprint.

\item[{\bf [BUK]}] D. M. Burns and S. G. Krantz, Rigidity of
holomorphic mappings and a new Schwarz lemma at the boundary,
{\it Jour. of the A.M.S.} 7(1994), 661-676.

\item[{\bf [CHE]}] D. Chelst, A generalized Schwarz lemma at the boundary,
{\it Proc.\ Amer.\ Math.\ Soc.} 129(2001), 3275--3278.

\item[{\bf [CIS1]}] J. A. Cima and T. J. Suffridge, Proper holomorphic
mappings from the two-ball to the three-ball, {\t Trans.\ Amer.\ Math.\
Soc.} 311 (1989), 227--239.

\item[{\bf [CIS2]}] J. A. Cima and T. J. Suffridge, Proper mappings between
balls in $\CC^n$, {\it Complex Analysis} (University Park, Pa., 1986),
66--82, Lecture Notes in Math., 1268, Springer, Berlin, 1987.

\item[{\bf [COH]}] R. Courant and D. Hilbert, {\it Methods of Mathematical
Physics}, Interscience, New York, 1953--1962.

\item[{\bf [DANG1]}] J. P. D'Angelo, The structure of proper rational
holomorphic maps between balls, {\it Several complex variables}
(Stockholm, 1987/1988), 227--244, Math. Notes, 38, Princeton Univ. Press,
Princeton, NJ, 1993.

\item[{\bf [DANG2]}] J. P. D'Angelo, The geometry of proper holomorphic
maps between balls, {\it The Madison Symposium on Complex Analysis}
(Madison, WI, 1991), 191--215, Contemp. Math., 137, Amer. Math. Soc.,
Providence, RI, 1992.

\item[{\bf [DANG3]}]  J. P. D'Angelo, Real hypersurfaces, orders of contact,
and applications, {\it Annals of Math.} 115(1982), 615-637.

\item[{\bf [DIF]}]  K. Diederich and J. E. Forn\ae ss, Pseudoconvex domains: 
Bounded strictly plurisubharmonic exhaustion functions, {\em Invent. Math.}
39(1977), 129-141.

\item[{\bf [FEF]}]  C. Fefferman, The Bergman kernel and biholomorphic mappings
of pseudoconvex domains, {\it Invent. Math.} 26(1974), 1-65.
			     
\item[{\bf [GRK]}]  R. E. Greene and S. G. Krantz, {\it Function Theory
of One Complex Variable}, $3^{\rm rd}$ ed., American Mathematical Society,
Providence, RI, 2006.
\end{enumerate}

\newpage

\begin{enumerate}
\item[{\bf [HUA]}] X. Huang, A boundary rigidity problem for holomorphic
mappings on some weakly pseudoconvex domains. {\it Canad.\ J. Math.}
47(1995), 405--420.

\item[{\bf [KAT]}]  Y. Katznelson, {\it An Introduction to Harmonic Analysis}, Dover, New York, 1976.

\item[{\bf [KRA1]}] S. G. Krantz, {\it Function Theory of Several Complex
Variables}, $2^{\rm nd}$ ed., American Mathematical Society, Providence,
RI, 2001.

\item[{\bf [KRA2]}] S. G. Krantz, {\it Cornerstones of Geometric Function
Theory: Explorations in Complex Analysis}, Birkh\"{a}user Publishing,
Boston, 2006.

\item[{\bf [KRA3]}] S. G. Krantz, The Kobayashi metric, extremal discs, \\
and biholomorphic mappings, preprint.

\item[{\bf [KRA4]}] S. G. Krantz, {\it Complex Analysis: The Geometric
Viewpoint}, A CARUS Monograph of the Mathematics Association of America,
Washington, D.C., 1990.
										 
\item[{\bf [KRA5]}] S. G. Krantz, A compactness principle in complex
analysis, {\it Division de Matematicas, Univ. Autonoma de Madrid
Seminarios}, vol. 3, 1987, 171-194.
										    
\item[{\bf [LEM]}]  L. Lempert, La metrique \mbox{K}obayashi et las representation des domains
sur la boule, {\em Bull. Soc. Math. France} 109(1981), 427-474.

\item[{\bf [LOW]}] K. L\"{o}wner, Untersuchungen \"{u}ber schlichte
konforme Abbildungen des Einheitskreises. I, {\it Math.\ Annalen}
89(1923), 103--121.

\item[{\bf [OSS]}] R. Osserman, A sharp Schwarz inequality on the boundary,
{\it Proc.\ Amer.\ Math.\ Soc.} 128(2000), 3513--3517.
		
\item[{\bf [SHO]}] D. Shoikhet, Another look at the
Burns-Krantz theorem, {\it J. Anal. Math.} 105(2008), 19--42.

\item[{\bf [VEL1]}]  J. Velling, Spherical Geometry and the Schwarzian Differential Equation,
thesis, Stanford University, 1985.

\item[{\bf [VEL2]}]  J. Velling, The uniformization of rectangles, an exercise in Schwarz's lemma,
{\it Amer.\ Math.\ Monthly} 39(1992), 112--115.

\item[{\bf [WEB]}] S. Webster, On mapping an $n$-ball into an $(n+1)$-ball
in complex space, {\it Pac. J. Math.} 81(1979), 267-272.
\end{enumerate}

\newpage

\begin{enumerate}
\item[{\bf [WU]}]  H. H. Wu, Normal families of holomorphic mappings, {\it Acta Math.} 119(1967),
193-233.

\item[{\bf [YAU]}]  S. T. Yau, A generalized Schwarz lemma for K\"{a}hler
manifolds, {\it Am. J. Math.} 100(1978), 197-204.

\end{enumerate}
\vspace*{.95in}

\noindent \begin{quote}
Department of Mathematics \\
Washington University in St. Louis \\
St.\ Louis, Missouri 63130 \\ 
{\tt sk@math.wustl.edu}
\end{quote}

\end{document}